\documentclass[reqno,12pt,a4paper,oneside]{article}
\usepackage{setspace}
\usepackage{amsfonts}
\usepackage{amsthm}
\usepackage[mathscr]{euscript}
\usepackage[cp1250]{inputenc}
\usepackage[T1]{fontenc}
\usepackage{amssymb}
\usepackage{amsmath}
\usepackage[dvipdf]{graphicx}
\def\={\discretionary{-}{-}{-}}

\newtheorem{tw}{Theorem}

\newtheorem{wn}{Corollary}

\title{Some inequalities for  Gurland's ratio of the gamma functions}
\author{Halina Wi\'sniewska}
\date{Institute of Mathematics,  Kazimierz Wielki University, 85-072 Bydgoszcz, Poland\\ email: halinkaw@ukw.edu.pl}
\begin{document}
\doublespacing
\maketitle

\vspace{12pt}

\begin{abstract}
This paper investigates the classical Gurland ratio of the gamma function and introduces its modified form, $\mathcal{G}^{\star}(x,y)$, which is particularly amenable to analytic expansions. By utilizing the Weierstrass product representation of the gamma function, we derive a finite expansion for the logarithm of $\mathcal{G}^{\star}(x,y)$ involving the Hurwitz zeta function. Explicit upper bounds for the remainder term are established, providing a rigorous basis for convergence analysis. As a direct consequence, we obtain new bilateral inequalities for the Gurland ratio and demonstrate the existence of a specific parameter $t(x,y)$ related to the Mean Value Theorem. Furthermore, we formulate open problems regarding the optimal localization of this parameter. These results extend the classical works of Gurland, Gautschi, and Merkle, offering new insights into the asymptotic behavior of gamma function ratios.
\end{abstract}


\section{Introduction} \label{wprowadzenie}

The Euler gamma function, $\Gamma(x)$, is undoubtedly one of the most important special functions in mathematical analysis, extending the factorial to complex numbers. A fundamental property of the gamma function is its logarithmic convexity for $x > 0$, a characteristic central to the Bohr-Mollerup theorem. This property implies various inequalities, which have been the subject of intensive research for decades.

In 1956, John Gurland \cite{Gurland} studied the ratio of gamma functions in the context of statistical distribution theory and derived the following inequality:
$$
\frac{\Gamma (x) \cdot \Gamma (y)}{{\Gamma}^2 (\frac{x+y}{2})} > 1
$$
for $x,y > 0$.

\noindent
The function on the left-hand side, commonly denoted as Gurland's ratio:
\begin{equation} \label{iloraz-Gurlanda}
\mathcal{G} (x,y) = \frac{\Gamma (x) \cdot \Gamma (y)}{{\Gamma}^2 (\frac{x+y}{2})}, \quad x, y > 0
\end{equation}
measures the deviation from the geometric mean of the gamma function values relative to the value at the arithmetic mean of the arguments. Since Gurland's original work, this ratio and its generalizations have been investigated by numerous mathematicians, including Gautschi \cite{Gautschi}, Kairies \cite{Kairies}, and Merkle \cite{Mercle}. For a comprehensive overview of these and other related inequalities for the gamma function, we refer the reader to the survey paper by Laforgia and Natalini \cite{Laforgia}. More recently, Chen and Choi \cite{Choi} extended these investigations by establishing the complete monotonicity of functions related to Gurland's ratio, paving the way for sharper bounds and deeper asymptotic analysis.

While inequalities provide bounds, exact expansions offer deeper insight into the asymptotic behavior of such ratios. In this paper, we aim to extend the classical results by providing a finite expansion for the logarithm of a modified Gurland ratio. To facilitate calculations and obtain a more symmetric form convenient for analytic expansions, we consider the modified ratio defined as:
\begin{equation}  \label{iloraz-Gurlanda-1}
{\mathcal{G}}_{\star} (x,y) = \frac{\Gamma (1+x) \cdot \Gamma (1+y)}{{\Gamma}^2 \left (1+ \frac{x+y}{2} \right )}.
\end{equation}

Using the Weierstrass product representation of the gamma function, we derive a finite expansion of $\ln \mathcal{G}_{\star}(x,y)$ involving the Hurwitz zeta function, $\zeta(s,a)$, whose definition and fundamental properties are discussed in detail in \cite{Wong}. We also provide error bounds for the remainder term and establish new bilateral inequalities that complement the work of Alzer \cite{Alzer} and Merkle \cite{Mercle}.

Let  $\zeta$ denote the Hurwitz zeta function defined as:
\begin{equation} \label{def-zeta}
 \zeta (a,x) = \sum_{n=0}^{\infty} {(n+x)}^{-a}.
 \end{equation}
 For the convenience of the reader and to facilitate a clear understanding of the presented results, all proofs are deferred to the final section of this paper.
\section{Main results}

We begin by establishing a finite expansion for the logarithm of the modified Gurland ratio, which provides a precise representation of the error term.

\begin{tw} \label{tw-zmi}
Let ${\mathcal{G}}_{\star} (x,y)$ denote the quotient given by $(\ref{iloraz-Gurlanda-1})$. Then for all $x, y > 0$  and integers $m \geq 2$:

\begin{equation} \label{dwie-gw}
\ln {\mathcal{G}}_{\star} (x,y) =  \sum_{k=1}^{m-1} \frac{1}{k} {\left (\frac{x-y}{2} \right )}^{2k} \zeta \left ( 2k, 1 + \frac{x+y}{2} \right ) + R_n^m (x,y),
\end{equation}
where $0 \leq R_n^m (x,y) \leq \varepsilon_n^m (x,y) = \frac{1}{m} {\left (\frac{x-y}{2} \right )}^{2m} \sum_{n=1}^{\infty} \frac{1}{{(n+x)}^m {(n+y)}^m}$.
\end{tw}
 
 Based on the error estimates derived in Theorem 1, we provide a sufficient condition for the expansion to converge to an infinite series representation.
 
\begin{wn} \label{wn1}
If $x, y >0$, then the condition $Q(x,y) = \frac{|x-y|}{2(1+ \sqrt{xy})} < 1$ implies that
$$
\ln {\mathcal{G}}_{\star} (x,y) =  S_{\infty} (x,y),
$$
where $S_{\infty} (x,y) = \lim_{m \to \infty} S_m (x,y)$ 
and 
$$
S_m (x,y) = \sum_{k=1}^{m-1} \frac{1}{k} {\left (\frac{x-y}{2} \right )}^{2k} \zeta \left ( 2k, 1 + \frac{x+y}{2} \right ),
$$
wherein
$$
0 \leq \ \ln {\mathcal{G}}_{\star} (x,y) - S_m (x,y) \leq \frac{1}{m} {\left ( \frac{|x-y|}{2(1+\sqrt{xy})} \right )}^{2m} \left (1+\frac{1+\sqrt{xy}}{2m-1} \right ).
$$

\end{wn}

Using the main theorem, we can also derive new bilateral inequalities that bound the logarithm of the ratio from above and below.
 
 \begin{wn} \label{wn2}
 Let $x, y > 0$. Then
 $$
{\left ( \frac{x-y}{2} \right )}^2 \cdot  \zeta (2, 1+\frac{x+y}{2}) \leq  \  \ln {\mathcal{G}}_{\star} (x,y) \ \leq {\left ( \frac{x-y}{2} \right )}^2 \cdot \zeta\left ( 2, 1+ \sqrt{xy} \right ).
 $$
 \end{wn}
 
 Furthermore, applying the intermediate value property allows us to express the ratio in terms of a specific parameter located between the geometric and arithmetic means.
 
 \begin{wn}  \label{wn3}
 For every $x, y > 0$ there exists $t := t(x,y) \in \left ( \sqrt{xy}, \frac{x+y}{2} \right )$ such that 
 $$
 \ln {\mathcal{G}}_{\star} (x,y) = {\left ( \frac{x-y}{2} \right )}^2 \cdot \zeta (2, 1+t).
 $$
 \end{wn}
 
 To conclude this section, we propose several questions that arise naturally from the convergence properties of the derived series.
 
 \noindent
 \textbf{Open problems:}
 
 \begin{itemize}
 \item[$1.$]  In Corollary \ref{wn1} we showed that $R_n^m (x,y)  \to 0$ under a sufficient condition. Does  $R_n^m (x,y)  \to 0$ hold for every $x,y > 0$?
 \item[$2.$] Does $\ln {\mathcal{G}}_{\star} (x,y) = S_{\infty} (x,y)$ hold generally, and if not, how large the difference  $\ln {\mathcal{G}}_{\star} (x,y) - S_{\infty} (x,y)$ can be?
 \item[$3.$] How is $t(x,y)$  located in the interval $\left (\sqrt{xy}, \frac{x+y}{2} \right )$ ? What is the distance from $t(x,y)$ to the midpoint of this interval?
 \end{itemize}
 
 \section{Proofs}
 
 \noindent
 \textbf{Proof of Theorem \ref{tw-zmi}}

We apply the following product representation (see \cite[p. 19]{Wong}):
\begin{equation} \label{gw1}
\Gamma (z) = \frac{1}{z} \prod_{n=1}^{\infty} {\left (  1+ \frac{1}{n} \right )}^z \cdot {\left (  1 + \frac{z}{n} \right ) }^{-1}, \ z \in \mathbb{C}.
\end{equation} 
Hence 

\begin{equation} \label{t2}
\Gamma (z+1) = z \Gamma (z) =  \prod_{n=1}^{\infty} {\left (  1+ \frac{1}{n} \right )}^z \cdot {\left (  1 + \frac{z}{n} \right ) }^{-1}, \ z \in \mathbb{C}.
\end{equation}
 
 Let $x,y > 0$.   By  ( \ref{iloraz-Gurlanda-1}),
we have 
 
\begin{equation} \label{p1}
{\mathcal{G}}_{\star} (x,y) = \frac{\Gamma (1+x) \cdot \Gamma (1+y)}{{\Gamma}^2 (1+\frac{x+y}{2})} = \frac{xy \Gamma (x) \Gamma (y)}{{\left ( \frac{x+y}{2} \right )}^2 {\Gamma}^2 \left ( \frac{x+y}{2} \right )}  = \frac{4xy}{{(x+y)}^2}{\mathcal{G}} (x,y).
\end{equation}

Now, let us consider the function $\mathcal{G}_{*}(x,y)$ using formula  (\ref{t2}). 
\begin{equation} \label{k2}
{\mathcal{G}}_{\star} (x,y) = \prod_{n=1}^{\infty} \frac{\frac{{\left ( 1+ \frac{1}{n} \right )}^x}{\left ( 1+ \frac{x}{n} \right )} \cdot \frac{{\left ( 1+ \frac{1}{n} \right )}^y}{\left ( 1+ \frac{y}{n} \right )}}{{\left ( \frac{{\left (1+ \frac{1}{n} \right )}^{\frac{x+y}{2}}}{1+ \frac{x+y}{2n}} \right )}^2} = \prod_{n=1}^{\infty} \frac{{\left (1+ \frac{x+y}{2n}  \right )}^2}{\left (  1+ \frac{x}{n} \right ) \cdot \left ( 1 + \frac{y}{n} \right )}.
\end{equation}

Set 
\begin{equation} \label{p2}
G_n (x,y) = \frac{{\left (1+ \frac{x+y}{2n}  \right )}^2}{\left (  1+ \frac{x}{n} \right ) \cdot \left ( 1 + \frac{y}{n} \right )}.
\end{equation}
Then
\begin{equation} \label{p3}
{\mathcal{G}}_{\star} (x,y) = \prod_{n=1}^{\infty} G_n (x,y).
\end{equation}

Considering the function $G_n (x,y)$:
$$
G_n (x,y) = \frac{{\left (1+ \frac{x+y}{2n}  \right )}^2}{\left (  1+ \frac{x}{n} \right ) \cdot \left ( 1 + \frac{y}{n} \right )} = {\left ( \frac{\left (  1+ \frac{x}{n} \right ) \cdot \left ( 1 + \frac{y}{n} \right )}{{\left (1+ \frac{x+y}{2n}  \right )}^2} \right ) }^{-1} =
$$
$$
 {\left ( \frac{1+ \frac{x+y}{n} + {\left (\frac{x+y}{2} \right )}^2 \frac{1}{n^2} + \frac{xy}{n^2} -{\left (\frac{x+y}{2} \right )}^2 \frac{1}{n^2}}{1+ \frac{x+y}{n} + {\left (\frac{x+y}{2} \right )}^2 \frac{1}{n^2}} \right )}^{-1} =
$$
\begin{equation} \label{k3}
{\left ( 1+ \frac{xy - {\left ( \frac{x+y}{2} \right )}^2}{n^2 + n(x+y) + {\left ( \frac{x+y}{2} \right )}^2 } \right )}^{-1} = {\left (1- \frac{{\left ( \frac{x-y}{2} \right )}^2}{{\left ( n+ \frac{x+y}{2} \right )}^2} \right )}^{-1}.
\end{equation}
Set 
\begin{equation} \label{k4}
c_n =  \frac{{\left ( \frac{x-y}{2} \right )}^2}{{\left ( n+ \frac{x+y}{2} \right )}^2}.
\end{equation}
By (\ref{k3}) we obtain
\begin{equation} \label{p4}
G_n (x,y) = {(1- c_n)}^{-1}.
\end{equation}

From (\ref{p4}) and  Maclaurin's theorem, we have
\begin{equation} \label{t5}
\ln G_n (x,y) = \ln {(1- c_n)}^{-1} = - \ln (1-c_n) = \sum_{k=1}^{m-1} \frac{c_n^k}{k} + \frac{1}{m} {(1- \theta_n)}^{-m} c_n^m,  \quad  \theta_n \in (0, c_n).
\end{equation}

Then, from (\ref{p3}) and (\ref{t5}):

$$
\ln {\mathcal{G}}_{\star} (x,y) =  \ln \left ( \prod_{n=1}^{\infty} G_n (x,y) \right ) = \sum_{n=1}^{\infty} \ln G_n (x,y) = \sum_{n=1}^{\infty}  \left ( \sum_{k=1}^{m-1} \frac{{c_n}^k}{k} + \frac{1}{m} {(1- \theta_n)}^{-m} c_n^m  \right )=
$$
\begin{equation} \label{p5}
\sum_{n=1}^{\infty} \sum_{k=1}^{m-1} \frac{{c_n}^k}{k} + \sum_{n=1}^{\infty} \frac{1}{m} {\left ( \frac{c_n}{1 - \theta_n} \right )}^m, \ \ \theta_n \in (0, c_n).
\end{equation}
Set
$$
S_n^m (x,y) =\sum_{n=1}^{\infty} \sum_{k=1}^{m-1} \frac{{c_n}^k}{k} 
$$
and
$$
R_n^m (x,y) = \sum_{n=1}^{\infty} \frac{1}{m} {\left ( \frac{c_n}{1 - \theta_n} \right )}^m.
$$
Then, for $ \theta_n \in (0, c_n)$ we have 
\begin{equation} \label{j-1}
\ln {\mathcal{G}}_{\star} (x,y) = S_n^m (x,y) + R_n^m (x,y).
\end{equation}

Since 
 \begin{equation} \label{p6}
 \frac{c_n}{1 - \theta_n}  < \frac{c_n}{1 - c_n} \stackrel{(\ref{k4})}{=} \frac{ \frac{{\left ( \frac{x-y}{2} \right )}^2}{{\left ( n+ \frac{x+y}{2} \right )}^2}}{1- \frac{{\left ( \frac{x-y}{2} \right )}^2}{{\left ( n+ \frac{x+y}{2} \right )}^2}} = \frac{{\left ( \frac{x-y}{2} \right )}^2}{(n+x)(n+y)},
 \end{equation}
 we obtain the bound:
 Now, let us consider the sum:

 $$
 R_n^m (x,y) = \sum_{n=1}^{\infty}\frac{1}{m} {\left ( \frac{c_n}{1 - \theta_n} \right )}^m \stackrel{(\ref{p6})}{<}  \sum_{n=1}^{\infty}\frac{1}{m} {\left ( \frac{{\left ( \frac{x-y}{2} \right )}^2}{(n+x)(n+y)} \right )}^m = 
 $$
 \begin{equation} \label{pp0}
  \frac{1}{m} {\left ( \frac{x-y}{2} \right )}^{2m}  \sum_{n=1}^{\infty} \frac{1}{{(n+x)}^m {(n+y)}^m}, \quad \theta_n \in (0, c_n).
 \end{equation}

 Changing the order of summation for $S_n^m (x,y)$:

$$
S_n^m (x,y) =\sum_{n=1}^{\infty} \sum_{k=1}^{m-1} \frac{{c_n}^k}{k} = \sum_{k=1}^{m-1} \sum_{n=1}^{\infty} \frac{{c_n}^k}{k} = \sum_{k=1}^{m-1} \left ( \frac{1}{k} \sum_{n=1}^{\infty} {c_n}^k \right ) \stackrel{(\ref{k4})}{=}
$$
$$
 \sum_{k=1}^{m-1} \left ( \frac{1}{k}  \sum_{n=1}^{\infty} \frac{{\left ( \frac{x-y}{2} \right )}^{2k}}{{\left ( n + \frac{x+y}{2} \right )}^{2k}} \right ) = \sum_{k=1}^{m-1} \left ( \frac{1}{k}  {\left ( \frac{x-y}{2} \right )}^{2k}  \sum_{n=1}^{\infty} \frac{1}{{\left ( n + \frac{x+y}{2} \right )}^{2k}} \right ) \stackrel{(\ref{def-zeta})}{=}
$$
 \begin{equation} \label{pp11}
\sum_{k=1}^{m-1}  \frac{1}{k}  {\left ( \frac{x-y}{2} \right )}^{2k}  \zeta \left ( 2k, 1+ \frac{x+y}{2} \right ).
 \end{equation}
 
 Combining (\ref{j-1}), (\ref{pp0}) and (\ref{pp11}) we obtain the assertion of our theorem:
 $$
 \ln {\mathcal{G}}_{\star} (x,y) = S_n^m (x,y) + R_n^m (x,y), 
 $$
 where 
 $$
 S_n^m (x,y) =\sum_{k=1}^{m-1}  \frac{1}{k}  {\left ( \frac{x-y}{2} \right )}^{2k}  \zeta \left ( 2k, 1+ \frac{x+y}{2} \right )
 $$
  and 
 
 $$
 0 \leq R_n^m (x,y) \leq  \varepsilon_n^m (x,y) = \frac{1}{m} {\left (\frac{x-y}{2} \right )}^{2m} \sum_{n=1}^{\infty} \frac{1}{{(n+x)}^m {(n+y)}^m}.
 $$
  As implied, by (\ref{p1}) 
 $$
 {\mathcal{G}}_{\star} (x,y) = \frac{4xy}{{(x+y)}^2} \mathcal{G} (x,y).
 $$ 
 then $\ln {\mathcal{G}} (x,y) = \ln  \left ( \frac{{(x+y)}^2}{4xy}  {\mathcal{G}}_{\star} (x,y) \right )$ and
 hence
 $$
 G(x,y) =  \frac{{(x+y)}^2}{4xy} \cdot e^{ S_n^m (x,y) +  R_n^m (x,y)},
 $$
 where $S_n^m (x,y)$ and $R_n^m(x,y)$ are given as above.
$\qed$

\noindent 
 \textbf{Proof of Corollary \ref{wn1}}
From Theorem \ref{tw-zmi},   for all $x, y > 0$ and integers $m \ge 2$, we have:

$$
\ln {\mathcal{G}}_{\star} (x,y) =  \sum_{k=1}^{m-1} \frac{1}{k} {\left (\frac{x-y}{2} \right )}^{2k} \zeta \left ( 2k, 1 + \frac{x+y}{2} \right ) + R_n^m (x,y),
$$
where the remainder satisfies:  
\begin{equation} \label{Rn}
0 \leq R_n^m (x,y) \leq \varepsilon_n^m (x,y) = \frac{1}{m} {\left (\frac{x-y}{2} \right )}^{2m} \sum_{n=1}^{\infty} \frac{1}{{(n+x)}^m {(n+y)}^m}.
\end{equation}
\vspace{11pt}

\noindent
Since $(n+x)(n+y) = n^2 + n(x+y) + xy$ and $x+y \ge 2 \sqrt{xy}$, we have
$$
(n+x)(n+y) \ge n^2 + 2n \sqrt{xy} + xy = {(n+ \sqrt{xy})}^2.
$$

 Applying this estimate to the sum, we obtain:
 $$
 \sum_{n=1}^{\infty} \frac{1}{{(n+x)}^m {(n+y)}^m} \leq   \sum_{n=1}^{\infty} \frac{1}{{(n+\sqrt{xy})}^{2m} } . 
 $$
 Using the integral test for convergence, we can bound the series as follows:
$$
\sum_{n=1}^{\infty} \frac{1}{(n+\sqrt{xy})^{2m}} \le \frac{1}{(1+\sqrt{xy})^{2m}} + \int_{1}^{\infty} \frac{ds}{(s+\sqrt{xy})^{2m}}=
$$
  \begin{equation} \label{pp1}
 \frac{1}{(1+\sqrt{xy})^{2m}} + \frac{1}{(2m-1)(1+\sqrt{xy})^{2m-1}}.
\end{equation}
  
  Consequently, using (\ref{Rn}) and (\ref{pp1}),  the remainder term is bounded by:
  
  $$
  R_{n}^{m}(x,y) \le \frac{1}{m} \left( \frac{x-y}{2} \right)^{2m} \left[ \frac{1}{(1+\sqrt{xy})^{2m}} + \frac{1}{(2m-1)(1+\sqrt{xy})^{2m-1}} \right].
  $$
  Let us denote this upper bound by $V_{m}(x,y)$. Simplifying the expression, we get:
 $$
 V_{m}(x,y) = \frac{1}{m} \left( \frac{x-y}{2(1+\sqrt{xy})} \right)^{2m} \left( 1 + \frac{1+\sqrt{xy}}{2m-1} \right).
 $$
 
 Hence
\begin{equation} \label{pp2}
R_n^m (x,y) \leq  V_n^m (x,y).
\end{equation}

If  we assume that $\frac{x-y}{2(1+ \sqrt{xy})} < 1$, then $V_m(x,y) \to 0$

Therefore,$$\lim_{m \to \infty} \left| \ln\mathcal{G}_{\star}(x,y) - S_{m}(x,y) \right| \le \lim_{m \to \infty} V_{m}(x,y) = 0,$$which implies $\ln\mathcal{G}_{\star}(x,y) = S_{\infty}(x,y)$. $\qed$

\noindent 
 \textbf{Proof of Corollary \ref{wn2}}
 
 From Theorem \ref{tw-zmi}, for every $x, y > 0$, we have the expansion: 
 $$
\ln {\mathcal{G}}_{\star} (x,y) =  \sum_{k=1}^{m-1} \frac{1}{k} {\left (\frac{x-y}{2} \right )}^{2k} \zeta \left ( 2k, 1 + \frac{x+y}{2} \right ) + R_n^m (x,y),
$$
where $0 \leq R_n^m (x,y) \leq \varepsilon_n^m (x,y) = \frac{1}{m} {\left (\frac{x-y}{2} \right )}^{2m} \sum_{n=1}^{\infty} \frac{1}{{(n+x)}^m {(n+y)}^m}$.
Since $R_{n}^{m}(x,y) \ge 0$, omitting the non-negative remainder yields the lower bound:
$$
\ln \mathcal{G}_{\star}(x,y) \ge \sum_{k=1}^{m-1} \frac{1}{k} \left( \frac{x-y}{2} \right)^{2k} \zeta\left( 2k, 1+\frac{x+y}{2} \right).
$$
Taking the first term of the sum (for $k=1$) gives the stated lower inequality:
$$
\ln \mathcal{G}_{\star}(x,y) \ge \left( \frac{x-y}{2} \right)^{2} \zeta\left( 2, 1+\frac{x+y}{2} \right).
$$
To prove the upper bound, we revisit the product representation from the proof of Theorem \ref{tw-zmi}:
$$
{\mathcal{G}}_{\star} (x,y) = \prod_{n=1}^{\infty} G_n (x,y),
$$
where 
$$
G_{n}(x,y) = \frac{(1+\frac{x+y}{2n})^{2}}{(1+\frac{x}{n})(1+\frac{y}{n})} = 1 + B_{n}(x,y),
$$
with
$$
B_{n}(x,y) = \frac{\frac{1}{n^2}(\frac{x-y}{2})^2}{1+\frac{x+y}{n}+\frac{xy}{n^2}} = \frac{(\frac{x-y}{2})^2}{n^2+n(x+y)+xy}.
$$
Using the inequality $\ln (1+u) \le u$ for $u > -1$, we have:
$$
\ln \mathcal{G}_{\star}(x,y) = \sum_{n=1}^{\infty} \ln (1 + B_{n}(x,y)) \le \sum_{n=1}^{\infty} B_{n}(x,y).
$$
Substituting the expression for $B_{n}(x,y)$:
$$
\sum_{n=1}^{\infty} B_{n}(x,y) = \left( \frac{x-y}{2} \right)^{2} \sum_{n=1}^{\infty} \frac{1}{n^2 + n(x+y) + xy}.
$$
Applying the inequality $x+y \ge 2\sqrt{xy}$ to the denominator:$$n^2 + n(x+y) + xy \ge n^2 + 2n\sqrt{xy} + xy = (n+\sqrt{xy})^2.$$Thus,
$$
\ln \mathcal{G}_{\star}(x,y) \le \left( \frac{x-y}{2} \right)^{2} \sum_{n=1}^{\infty} \frac{1}{(n+\sqrt{xy})^2} = \left( \frac{x-y}{2} \right)^{2} \zeta(2, 1+\sqrt{xy}).
$$
This completes the proof. 
$\qed$

\noindent 
 \textbf{Proof of Corollary \ref{wn3}}
 Let us define the function $f(t)$ for $t > 0$ as:
 $$
 f(t) = \left( \frac{x-y}{2} \right)^{2} \zeta(2, 1+t).
 $$
 The function $t \mapsto \zeta(2, 1+t)$ is continuous and strictly decreasing for $t > 0$.From Corollary  \ref{wn2}, we have the inequalities:$$f\left( \frac{x+y}{2} \right) \le \ln \mathcal{G}_{\star}(x,y) \le f(\sqrt{xy}).
 $$
 Since $f(t)$ is continuous on the interval $[\sqrt{xy}, \frac{x+y}{2}]$, by the Intermediate Value Theorem, there exists a value $t(x,y) \in [\sqrt{xy}, \frac{x+y}{2}]$ such that:
 $$
 \ln \mathcal{G}_{\star}(x,y) = f(t(x,y)).$$This concludes the proof.  $\qed$

\noindent
\textbf{Acknowledgement}.   The author thanks to Professor Marek W\'ojtowicz for remarks and comments which improved the quality of this paper.

\end{document}